# On the average uncertainty for systems with nonlinear coupling

*Kenric P. Nelson[1,2], Sabir R. Umarov[3], and Mark A. Kon[1]*

## Abstract

The increased uncertainty and complexity of nonlinear systems have motivated investigators to consider generalized approaches to defining an entropy function. New insights are achieved by defining the average uncertainty in the probability domain as a transformation of entropy functions. The Shannon entropy when transformed to the probability domain is the weighted geometric mean of the probabilities. For the exponential and Gaussian distributions, we show that the weighted geometric mean of the distribution is equal to the density of the distribution at the location plus the scale (i.e. at the width of the distribution). The average uncertainty is generalized via the weighted generalized mean, in which the moment is a function of the nonlinear source. Both the Rényi and Tsallis entropies transform to this definition of the generalized average uncertainty in the probability domain. For the generalized Pareto and Student's t-distributions, which are the maximum entropy distributions for these generalized entropies, the appropriate weighted generalized mean also equals the density of the distribution at the location plus scale. A coupled entropy function is proposed, closely related to the normalized Tsallis entropy, but incorporating a distinction between the additive coupling and multiplicative coupling.

## 1 Introduction

Entropy is the standard-bearer for defining the average uncertainty of a probability distribution or density function [1–3]. Boltzmann, Gibbs, and Shannon (BGS) demonstrated that the kernel $-\ln p$ is necessary to form a weighted arithmetic average of the uncertainty of a probability distribution. Khinchin showed that the entropy function $H = -\sum_{i=1}^{N} p_i \ln p_i$ is unique given the axioms that it is

a) continuous in $p_i$,
b) maximized at the uniform distribution,
c) not changed by a state with zero probability, and
d) additive for conditionally independent probabilities.

Correspondence: Kenric Nelson, kenricpn@bu.edu, 781-645-8564
[1] Boston University, Boston, MA
[2] Raytheon Company, Woburn, MA
[3] University of New Haven, New Haven, CT

Rényi, Tsallis and Amari [4–6] sought to broaden the definition of average uncertainty to account for the influence of nonlinear dynamics in complex systems. Generalized entropy measures have been applied to a variety of nonlinear systems such as decision making under risk [7,8], communication channel equalization [9], compressive sensing [10], the edge of chaos [11], space plasma [12,13], high energy physics [14,15] and quantum entanglement [16]. Hanel and Thurner [17,18] have shown that requiring just the first three of Khinchin's axioms leads to a two-parameter generalization of entropy, one of which is the Tsallis generalization and the focus of this communication.

The objective of the paper is to show the generalized average uncertainty for a nonlinear system can be defined in the probability domain as the weighted generalized mean. This is derived as a transformation of the generalized entropy functions. New evidence for the importance of generalizing the average uncertainty is provided. For the distributions that minimize the weighted generalized mean (maximum generalized entropy) constrained by a location and scale, the density at the location plus the scale is the generalized average uncertainty.

A review of the average uncertainty for important members of the exponential family provides a helpful framework prior to introduction of the effects of nonlinearity. The Boltzmann-Gibbs-Shannon entropy transformed to the probability domain, is the weighted geometric mean of the distribution, that is

$$P_{avg} \equiv \exp\left(+\sum_{i=1}^{N} p_i \ln p_i\right) = \prod_{i=1}^{N} p_i^{p_i}, \quad (1)$$

where the weights, now as powers, are also probabilities. The *average uncertainty* $P_{avg}$ is the average probability of being able to determine the state of the system. The average uncertainty ranges from certainty ($P_{avg}=1$) when one of the states is certain, to ($P_{avg}=1/N$) when all states are equiprobable. The intuition is that the total probability of the distribution is the product of the probabilities and the average is determined by weighting each term by the probability. Assuming a continuous distribution, the average density for the exponential distribution is then

$$f_{avg}\left(\tfrac{1}{\sigma}e^{\frac{-(x-\mu)}{\sigma}}\right) = \exp\left(\int_{\mu}^{\infty} \tfrac{1}{\sigma}e^{\frac{-(x-\mu)}{\sigma}} \ln\left(\tfrac{1}{\sigma}e^{\frac{-(x-\mu)}{\sigma}}\right) dx\right)$$
$$= \frac{1}{\sigma e} \quad (2)$$

and for the Gaussian distribution the average density is

$$f_{avg}\left(\tfrac{1}{\sqrt{2\pi}\sigma}e^{-\tfrac{1}{2}\left(\frac{x-\mu}{\sigma}\right)^2}\right) = \exp\left(\int_{-\infty}^{\infty} \tfrac{1}{\sqrt{2\pi}\sigma}e^{-\tfrac{1}{2}\left(\frac{x-\mu}{\sigma}\right)^2} \ln\left(\tfrac{1}{\sqrt{2\pi}\sigma}e^{-\tfrac{1}{2}\left(\frac{x-\mu}{\sigma}\right)^2}\right) dx\right)$$
$$= \frac{1}{\sqrt{2\pi e}\sigma}. \quad (3)$$

In both cases the average density is equal to $f(\mu+\sigma)$.

Rényi showed that information theory could be broadened by considering the weighted generalized mean of probabilities [19] as the kernel prior to applying the logarithm function. Through consideration of the statistics of a weighted distribution $p_i^q$, Tsallis showed that use of both the weighted generalized mean and a generalization of the logarithm function provided a model of non-additive entropy. Because of its role in raising probabilities to a power, the parameter *q* can be interpreted as the number of independent random variables whose combined distribution provides a basis for determining generalized statistical properties. This analysis, broadly referred to as nonextensive statistical mechanics [20], has been shown to model a wide variety of complex systems.

Section 2 introduces the concept of nonlinear statistical coupling which is an interpretation of nonextensive statistical mechanics focused on the role of nonlinearity in deforming statistical analysis. From examination of multivariate distributions, it was shown that the coupling parameter $\kappa$ is related to *q* via the influence of the dimensions *d* and power $\alpha$ of the state variable by $q = 1 + \frac{\alpha\kappa}{1+d\kappa}$. Other approaches have been taken to parameterizing the Tsallis statistics, including $\kappa' = 1-q$, which was originally proposed as a definition of nonlinear statistical coupling [21] and has been utilized by other investigators [22,23]; and $\kappa'' = \frac{1}{1-q}$ which is the translation to the kappa-distribution used in space physics [12]. There are also generalizations of entropy, which use the kappa symbol, such as the work by Kaniadakis[24].

In Section 3, the main objective relating average uncertainty and the width of coupled exponential distributions is established. In Section 4 we show that the generalized entropies of Renyi and Tsallis along with a modified normalization of the Tsallis entropy denoted as the coupled entropy, can each be expressed as the weighted generalized mean of the distribution and a transformation to the entropy scale using a generalized logarithm function.

## 2   Nonlinear Statistical Coupling

Fundamental to complex systems is the influence of nonlinearity on dynamics [25,26]. As such in modeling the statistical mechanics of a complex system, we choose the source of nonlinearity to define the *nonlinear statistical coupling (NSC)* $\kappa$ [21]. NSC has two contexts, external coupling between random variables and internal coupling between the states of a random variable [27]. The role of the coupled product function in modeling external coupling was developed in [28], building upon the q-algebra of nonextensive statistical mechanics [29]. In this paper, the focus is on internal coupling. The probability average translated from BGS entropy in (1) treats the probability of each state as independent, thus multiplying the probabilities of the distribution. The dependence due to the probabilities summing to one is incorporated via the power term. The generalized average uncertainty introduced in Section 3 will make use of the coupled product as a model of dependency between the states. First we introduce the basic relationships of nonlinear statistical coupling.

Positive coupling $(0<\kappa<\infty)$ is the domain of heavy-tail statistics in which the nonlinearity causes increased variation. Examples include multiplicative noise [30], superstatistics in which the variance fluctuates [31,32], and the inverse of the degree of freedom [33–35]. Negative coupling $(-1<\kappa<0)$ is the domain of compact-support statistics in which the nonlinearity causes reduction in variation. The two domains are related by a conjugate dual $\frac{-\kappa}{1+\kappa}$. Linear systems $(\kappa=0)$ are the domain of exponential statistics and logarithmic measures of information.

The fluctuation of the variance, studied as superstatistics [31,32], is a helpful example. For a system governed by an exponential distribution,

$$\frac{1}{\sigma'}e^{-\frac{x}{\sigma'}}, \; x\geq 0, \tag{4}$$

but with the inverse scale $\beta'=\frac{1}{\sigma'}$ fluctuating according to a gamma distribution the modified statistics are governed by a deformed exponential distribution, which we will call a coupled-exponential. The scale and shape of the coupled exponential are defined by the mean and relative variance of the inverse scale

$$\sigma = \langle \beta' \rangle^{-1},$$
$$\kappa = \frac{\langle \beta'^2 \rangle - \langle \beta' \rangle^2}{\langle \beta' \rangle^2}. \tag{5}$$

For this mean and relative variance, the gamma distribution is

$$f(\beta')=\frac{1}{\Gamma\left(\frac{1}{\kappa}\right)\left(\frac{\kappa}{\sigma}\right)}(\beta')^{\frac{1}{\kappa}-1}e^{-\frac{\sigma\beta'}{\kappa}}, \tag{6}$$

and the coupled exponential distribution is

$$\frac{1}{\sigma}(1+\kappa x)^{-\left(\frac{1}{\kappa}+1\right)}=\int_0^\infty f(\beta')(\beta' e^{-\beta' x})d\beta'. \tag{7}$$

The *coupled exponential* function is notated as an analog of the exponential function

$$\exp_\kappa(x)\equiv(1+\kappa x)_+^{\frac{1}{\kappa}+1}, \; (a)_+ \equiv \max(0,a). \tag{8}$$

The inverse of the exponent is a modified coupling, which will be referred to as the *multiplicative coupling*. For $d$ dimensions [28,36] the exponent generalizes to $\kappa_d^{-1}\equiv\frac{1}{\kappa}+d$, but only the one-dimensional function will be used in this paper. $\kappa_d$ will be used to modify the product function. The inverse of the multiplicative coupling is equal to the asymptotic fractal dimension $\lim_{x\to\infty}(1+\kappa x)^{-\kappa_d^{-1}}=\mathcal{O}\left(x^{-\kappa_d^{-1}}\right)$ of the coupled exponential distribution (7); thus the inverse of the coupling terms are summarized by

$$\kappa^{-1} = \nu \quad \text{Degree of Freedom,}$$
$$\kappa_d^{-1} = D \quad \text{Asymptotic Fractal Dimension.} \tag{9}$$

Raising to a power is not identical to multiplying the argument, so an additional parameter definition is required to show the relationships. In anticipation of its use with distributions, the power used is $1/\alpha$

$$\exp_\kappa^{1/\alpha}(x) \equiv (1+\kappa x)_+^{\frac{1}{\alpha}\left(\frac{1}{\kappa}+1\right)},$$
$$\exp_{\alpha,\kappa}(x) \equiv (1+\alpha\kappa x)_+^{\frac{1}{\alpha}\left(\frac{1}{\kappa}+1\right)}. \tag{10}$$

These notational definitions are needed to clarify that raising to a power is equivalent to both multiplying the argument and dividing the nonlinear coupling terms

$$\exp_\kappa^{1/\alpha}(x) = \exp_{\alpha,\kappa}(x/\alpha). \tag{11}$$

The nonlinear statistical coupling has the effect of changing the relationship between variables in the state space. Multiplication of independent coupled exponential distributions results in a combined variable of

$$(1+\kappa x)(1+\kappa y) = (1+\kappa(x+y+\kappa xy)). \tag{12}$$

The nonlinear term is viewed as forming the *coupled addition* [21,29]

$$x \oplus_\kappa y = x + y + \kappa xy, \tag{13}$$

which via the *coupled subtraction* is seen as a dilation of the state space

$$x \ominus_\kappa y \equiv \frac{x-y}{1+\kappa y}. \tag{14}$$

The inverse of the coupled exponential for the domain $1+\kappa x > 0$ is the *coupled logarithm*

$$\ln_\kappa(x) \equiv \frac{1}{\kappa}\left(x^{\frac{\kappa}{1+\kappa}} - 1\right), \quad x > 0. \tag{15}$$

Importantly, the integral over the unit interval is invariant to the deformation $\int_0^1 \ln_\kappa(x^{-1})dx = 1$. This insures that the coupled logarithm deforms the measure of information without modifying the "total information" over the unit interval. Incorporating, the power $\alpha$ the definition is

$$\ln_{\alpha,\kappa}(x) \equiv \frac{1}{\alpha}\ln_\kappa(x^\alpha) \equiv \frac{1}{\alpha\kappa}\left(x^{\frac{\alpha\kappa}{1+\kappa}} - 1\right), \quad x > 0 \tag{16}$$

The role of the coupled addition in modifying the state space can now be expressed in terms of the coupled exponential and logarithm functions

$$\prod_{i=1}^{N} \exp_\kappa(x_i) = \exp\left(\sum_{i=1}^{N} {}_{\oplus_\kappa} x_i\right), \quad \sum_{i=1}^{N} {}_{\oplus_\kappa} x_i \equiv x_1 \oplus_\kappa \ldots \oplus_\kappa x_N,$$
$$\ln_\kappa\left(\prod_{i=1}^{N} x_i\right) = \sum_{i=1}^{N} {}_{\oplus_\kappa} \ln_\kappa x_i.$$
(17)

The coupled product has a complimentary role, but operates on the exponent of the variable. Thus not only do the coupled sum and product not form an algebra with a distribution property [29], but are more accurately two complementary algebras [22,37]. The *coupled product* of positive-valued variables $f > 0$ and $g > 0$ is

$f \otimes_{\kappa_1} g \equiv \left(f^{\kappa_1} + g^{\kappa_1} - 1\right)_+^{\frac{1}{\kappa_1}}$, where the multiplicative coupling $\kappa_1 = \frac{\kappa}{1+\kappa}$ is the inverse of the exponent of the coupled exponential function. This definition extends to multiple variables and has the following properties

$$\prod_{i=1}^{N} {}_{\otimes_{\frac{\kappa}{1+\kappa}}} \exp_\kappa(x_i) \equiv \left(\sum_{i=1}^{N} \exp_\kappa^{\frac{\kappa}{1+\kappa}}(x_i) - (N-1)\right)^{\frac{1+\kappa}{\kappa}}$$
$$= \exp_\kappa\left(\sum_{i=1}^{N} x_i\right), \text{ for } \exp_\kappa(x_i) > 0 \quad (18)$$
$$\ln_\kappa\left(\prod_{i=1}^{N} {}_{\otimes_{\frac{\kappa}{1+\kappa}}} x_i\right) = \sum_{i=1}^{N} \ln_\kappa x_i, \text{ for } x_i > 0.$$

In [28] the formation of multivariate coupled exponentials required the coupled product to be defined such that the output dimension was the sum of the input dimensions. This paper is focused on the coupling between states for a single random variable and as such the dimension does not change. The *coupled power* will also be needed,

$$x^{\otimes_{\kappa_1}^a} \equiv \left(ax^{\kappa_1} - (a-1)\right)_+^{\frac{1}{\kappa_1}}, \; x > 0. \quad (19)$$

Generalization of the stretched exponential distributions $f(x) \propto \exp\left(\frac{-x^\alpha}{\alpha \sigma^\alpha}\right)$ utilizes the two-parameter definitions from Equation (10). Thus the coupled stretched exponential is expressed as

$$\exp_\kappa^{-1/\alpha}\left(\frac{x^\alpha}{\sigma^\alpha}\right) = \exp_{-\alpha,\kappa}\left(\frac{-x^\alpha}{\alpha \sigma^\alpha}\right) = \left(1 + \kappa \frac{x^\alpha}{\sigma^\alpha}\right)_+^{\frac{1+\kappa}{-\alpha\kappa}}. \quad (20)$$

The coupled exponential $(\alpha = 1)$ and coupled Gaussian $(\alpha = 2)$ distributions are of particular interest given their role as maximum generalized entropy distributions with constraints on the scale of the distribution. These distributions [38] are defined as

$$\frac{1}{Z(\kappa,\sigma,\alpha)}\exp_\kappa^{-1/\alpha}\left(\left(\frac{x-\mu}{\sigma}\right)^\alpha\right), \quad \kappa \geq -1. \tag{21}$$

For $\kappa \geq 0$ and $\alpha = 1$ the domain of the random variable is $x \geq \mu$ and $Z(\kappa,\sigma,\alpha) = \sigma$.

For $\kappa \geq 0$ and $\alpha = 2$, $-\infty < x < \infty$ and $Z(\kappa,\sigma,\alpha) = \dfrac{\sqrt{\kappa}\left(\frac{1+\kappa}{2\kappa}\right)!}{\sqrt{\pi}\sigma(1+\kappa)\left(\frac{1}{2\kappa}\right)!}$, where

$x! \equiv \Gamma(x+1)$ is Euler's gamma function. The parameters of the distributions are the location $\mu$, the scale $\sigma$, and the normalization $Z$. One of the benefits of this approach to the definition is that the source of coupling $\kappa$ corresponds exactly to the shape parameter for the generalized Pareto distribution and the inverse of the degree of freedom for the Student's t-distribution. For the compact-support domain, the degree of freedom is negative and less than -1.

**Definition 1 Coupled Probability/Density** Given a discrete probability distribution $\mathbf{p} = \{p_i, i = 1,...,N\}$ the following distribution is called the coupled probability distribution

$$P_i^{(\alpha,\kappa)} = \frac{p_i^{1+\frac{\alpha\kappa}{1+\kappa}}}{\sum_{j=1}^N p_j^{1+\frac{\alpha\kappa}{1+\kappa}}} = \frac{p_i \Big/ \prod_{\substack{j=1 \\ j \neq i}}^N p_j^{\frac{\alpha\kappa}{1+\kappa}}}{\sum_{i=1}^N \left( p_i \Big/ \prod_{\substack{j=1 \\ j \neq i}}^N p_j^{\frac{\alpha\kappa}{1+\kappa}} \right)}. \tag{22}$$

Similarly, for a given continuous density $f(x)$ the coupled density function is defined by

$$f^{(\alpha,\kappa)}(x) \equiv \frac{f^{1+\frac{\alpha\kappa}{1+\kappa}}(x)}{\int_{-\infty}^{\infty} f^{1+\frac{\alpha\kappa}{1+\kappa}}(x)dx}. \tag{23}$$

The integral in (23) is assumed finite. The expression on the right of Equation (22) is shown to make evident the modeling of a coupled state of the system. The Cauchy distribution $(\alpha = 2, \kappa = 1)$ is an illustrative example. In this case $\frac{\alpha\kappa}{1+\kappa} = 1$ and the coupled-probability is formed by dividing the probability $p_i$ of state $i$ by the product of all the probabilities of all the other states and then renormalizing them. This procedure has a direct connection with the origin of the Cauchy distribution as the division of two random variables.

An important result of nonextensive statistical mechanics is the relationship between coupled-moments and the parameters of the coupled-exponential and coupled-Gaussian distributions [21,39]. Using the coupled algebra the generalization of the $n^{th}$-moment takes the form

$$\left\langle x_i^n \right\rangle_\kappa \equiv \sum_{i=1}^N x_i^n P_i^{(n,\kappa)} = \frac{\sum_{i=1}^N x_i^n p_i^{1+n\frac{\kappa}{1+\kappa}}}{\sum_{i=1}^N p_i^{1+n\frac{\kappa}{1+\kappa}}}. \tag{24}$$

The $\alpha$ parameter is not required because it drops out from the expression $\frac{n}{\alpha}\frac{\alpha\kappa}{1+\kappa}$.

**Lemma 1** Given either the coupled-exponential or the coupled-Gaussian distribution, the generalized scale parameter $\sigma$ is determined by the generalizations of the 1st and 2nd moments, respectively:

$$\left\langle x \right\rangle_\kappa = \int_{-\infty}^{\infty} x \left( \frac{1}{\sigma} \exp_\kappa^{-1}\left(\frac{x}{\sigma}\right) \right)^{\frac{1+2\kappa}{1+\kappa}} dx \bigg/ \int_{-\infty}^{\infty} \left( \frac{1}{\sigma} \exp_\kappa^{-1}\left(\frac{x}{\sigma}\right) \right)^{\frac{1+2\kappa}{1+\kappa}} dx = \sigma, \tag{25}$$

$$\left\langle x^2 \right\rangle_\kappa = \int_{-\infty}^{\infty} x^2 \left( \frac{1}{Z} \exp_\kappa^{-1/2}\left(\frac{x^2}{\sigma^2}\right) \right)^{\frac{1+3\kappa}{1+\kappa}} dx \bigg/ \int_{-\infty}^{\infty} \left( \frac{1}{Z} \exp_\kappa^{-1/2}\left(\frac{x^2}{\sigma^2}\right) \right)^{\frac{1+3\kappa}{1+\kappa}} dx = \sigma^2. \tag{26}$$

See [39] for the proof and [21] for the notation.

## 3   The Coupled Average Uncertainty

The coupled algebra provides a foundation for defining a generalization of the average uncertainty given coupling of statistical states which is equal to the weighted generalized mean (or weighted p-norm) [22,27,40]. Equation (1) is generalized using the coupled exponential and coupled logarithm functions.

**Definition 2 Coupled Log Average**  Given a set of weights $\{w_i,\ i=1,..N\}$ such that $\sum_{i=1}^N w_i = 1$ and a set of variables $\{x_i,\ i=1,..N\}$, the coupled-log average of the variables is

$$\exp_\kappa^{-1/\alpha}\left\{\sum_{i=1}^N w_i \left(\ln_\kappa x_i^{-\alpha}\right)\right\}. \tag{27}$$

**Lemma 2** The coupled-log average

   a) can be expressed equivalently using the coupled product, and

   b) is equal to the weighted generalized mean (or weighted p-norm).

Proof a) We have

$$\exp_\kappa^{-1/\alpha}\left\{\sum_{i=1}^{N} w_i\left(\ln_\kappa x_i^{-\alpha}\right)\right\} = \prod_{i=1}^{N} {}_{\otimes_{\frac{-\alpha\kappa}{1+\kappa}}} x_i^{\otimes_{\frac{-\alpha\kappa}{1+\kappa}}^{w_i}}, \qquad (28)$$

where we used the properties in Equations (18) and (19).

b) We have

$$\exp_\kappa^{-1/\alpha}\left\{\sum_{i=1}^{N} w_i\left(\ln_\kappa x_i^{-\alpha}\right)\right\}$$

$$= \left(1+\kappa\left(\sum_{i=1}^{N}\frac{w_i}{\kappa}\left(x_i^{\frac{-\alpha\kappa}{1+\kappa}}-1\right)\right)\right)^{\frac{1+\kappa}{-\alpha\kappa}} \qquad (29)$$

$$= \left(\sum_{i=1}^{N} w_i x_i^{\frac{-\alpha\kappa}{1+\kappa}}\right)^{\frac{1+\kappa}{-\alpha\kappa}}.$$

Given the central role of the parameter $m = \dfrac{\alpha\kappa}{1+\kappa}$ for the weighted generalized mean, the coupled log average is summarized as

$$\exp_\kappa^{-1/\alpha}\left\{\sum_{i=1}^{N} w_i\left(\ln_\kappa x_i^{-\alpha}\right)\right\} = \left(\sum_{i=1}^{N} w_i x_i^{-m}\right)^{\frac{-1}{m}}, \quad m = \frac{\alpha\kappa}{1+\kappa}. \qquad (30)$$

The coupled-log average extends to continuous variables by the relationship

$$\exp_\kappa^{-1/\alpha}\left\{\int_{-\infty}^{\infty} w(x)\ln_\kappa f^{-\alpha}(x)\,dx\right\} = \left(\int_{-\infty}^{\infty} w(x) f^{-m}(x)\,dx\right)^{\frac{-1}{m}}. \qquad (31)$$

Using the coupled probability as the weight of the generalized mean we can now define the average uncertainty of a distribution that originates from a system with nonlinear coupling.

**Definition 3 Coupled Average Uncertainty** Given a distribution $\mathbf{p} = \{p_i,\ i=1,\ldots,N\}$ its coupled average uncertainty is defined by the coupled log average with $x_i = p_i$ and $w_i = P_i^{(\alpha,\kappa)}$

$$P_{\kappa avg} \equiv \left(\sum_{i=1}^{N} P_i^{(\alpha,\kappa)} p_i^{\frac{-\alpha\kappa}{1+\kappa}}\right)^{\frac{1+\kappa}{-\alpha\kappa}}. \qquad (32)$$

Furthermore, the coupled average uncertainty simplifies to the weighted generalized mean with the distribution as the weight; thus

$$P_{\kappa avg} \equiv \left( \frac{\sum_{i=1}^{N} p_i^{1+m} p_i^{-m}}{\sum_{j=1}^{N} p_j^{1+m}} \right)^{\frac{-1}{m}} = \left( \sum_{i=1}^{N} p_i^{1+m} \right)^{\frac{1}{m}}, \ m = \frac{\alpha \kappa}{1+\kappa}, \quad (33)$$

and for a continuous distribution

$$f_{\kappa avg} = \exp_\kappa^{-1/\alpha} \left\{ \frac{\int_{-\infty}^{\infty} f^{1+\frac{\alpha\kappa}{1+\kappa}}(x) \ln_\kappa f^{-\alpha}(x) dx}{\int_{-\infty}^{\infty} f^{1+\frac{\alpha\kappa}{1+\kappa}}(x) dx} \right\} = \left( \int_{-\infty}^{\infty} f^{1+m}(x) dx \right)^{\frac{1}{m}}. \quad (34)$$

Applying the coupled average uncertainty to the coupled exponential and coupled Gaussian distributions, shown in Figure 1, produces the main result regarding its relationship with the location and scale of the distribution.

**Theorem 1** The coupled average $f_{\kappa avg}$ of either the coupled exponential distribution or the coupled Gaussian distribution with a weight defined by the coupled probability is equal to the density of the distribution evaluated at $x = \mu + \sigma$.

Proof: Without loss of generality we set $\mu = 0$. For the coupled exponential distribution, which corresponds to $\alpha = 1$ in Equation (21), that is $f(x) = \frac{1}{\sigma} \exp_\kappa^{-1} \left( \frac{x}{\sigma} \right)$, the coupled average uncertainty is

$$\begin{aligned} f_{\kappa avg} &= \left( \int_0^\infty f^{1+m}(x) dx \right)^{\frac{1}{m}} \\ &= \sigma^{\frac{1+2\kappa}{\kappa}} \left( \int_0^\infty \left( 1 + \kappa \frac{x}{\sigma} \right)^{\frac{1+2\kappa}{-\kappa}} dx \right)^{\frac{1+\kappa}{\kappa}} \\ &= \sigma^{\frac{1+2\kappa}{\kappa}} \left( (1+\kappa)^{-1} \sigma \right)^{\frac{1+\kappa}{\kappa}} \\ &= \frac{1}{\sigma} e_\kappa^{-1}(1) = f(\mu + \sigma). \end{aligned} \quad (35)$$

Similarly, for the coupled Gaussian distribution, which corresponds to $\alpha=2$ in Equation (21), that is $f(x)=\dfrac{1}{Z(\kappa,\sigma,2)}\exp_\kappa^{-1/2}\left(\dfrac{x^2}{\sigma^2}\right)$, the coupled average uncertainty is

$$\begin{aligned}
f_{\kappa avg} &= \left(\int_{-\infty}^{\infty} f^{1+m}(x)dx\right)^{\frac{1}{m}} \\
&= \left(\frac{\sqrt{\kappa}\left(\frac{1+\kappa}{2\kappa}\right)!}{\sqrt{\pi}\sigma(1+\kappa)\left(\frac{1}{2\kappa}\right)!}\right)^{\frac{1+3\kappa}{2\kappa}} \left(\int_{-\infty}^{\infty}\left(1+\kappa\frac{x^2}{\sigma^2}\right)^{\frac{1+3\kappa}{-2\kappa}}dx\right)^{\frac{1+\kappa}{2\kappa}} \\
&= \left(\frac{\sqrt{\kappa}\left(\frac{1+\kappa}{2\kappa}\right)!}{\sqrt{\pi}\sigma(1+\kappa)\left(\frac{1}{2\kappa}\right)!}\right)^{\frac{1+3\kappa}{2\kappa}} \left(\frac{\sqrt{\pi}\sigma\left(\frac{1}{2\kappa}\right)!}{\sqrt{\kappa}\left(\frac{1+\kappa}{2\kappa}\right)!}\right)^{\frac{1+\kappa}{2\kappa}} \quad (36) \\
&= \left(\frac{\sqrt{\kappa}\left(\frac{1+\kappa}{2\kappa}\right)!}{\sqrt{\pi}\sigma(1+\kappa)\left(\frac{1}{2\kappa}\right)!}\right)(1+\kappa)^{-\frac{1+\kappa}{2\kappa}} \\
&= \frac{1}{Z(\kappa,\sigma,2)}e_\kappa^{-1/2}(1) = f(\mu+\sigma)
\end{aligned}$$

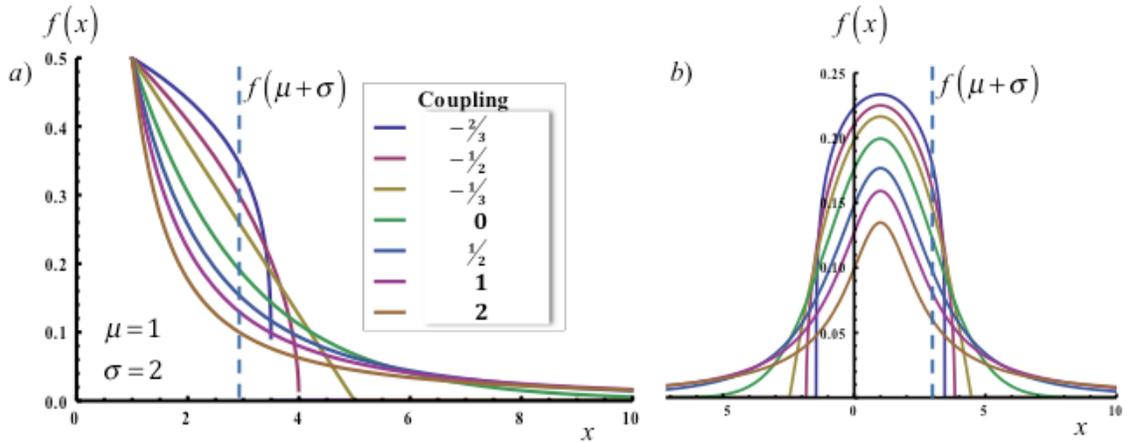

**Figure 1** The a) coupled exponential and b) coupled Gaussian distributions have the property that the coupled average uncertainty, equation (34), is equal to the density at the location plus the scale $f(\mu+\sigma)$. Several examples of the distributions with coupling ranging from $-\tfrac{2}{3}$ to 2 are shown.

Theorem 1 has important implications for defining the generalized average uncertainty. For $\kappa \to 0$ the weighted generalized mean is the weighted geometric mean. Just as the standard deviation is an expression of the variation of the random variable, the geometric mean of a distribution weighted by the distribution is an expression of average uncertainty. For the Gaussian distribution, one can verify that the weighted geometric mean of the distribution is equal to the density evaluated at $x = \mu + \sigma$.

The average uncertainty is generalized by using the coupled probability distribution for the weight and taking the $-\alpha\kappa/(1+\kappa)$ moment. And when $f(x)$ is a member of the coupled distributions $f(\mu+\sigma) = f_{\kappa avg}$. This relationship is shown in

Figure 2 for the heavy-tail domain of the coupled-Gaussian distribution. The coupled average uncertainty versus the source of coupling parameter for the distribution is plotted over the range $0 \leq \kappa \leq 1$, which includes the Gaussian $(\kappa=0)$, the boundary with infinite variance $(\kappa=0.5)$, and the Cauchy distribution $(\kappa=1)$. Four values for the coupled average uncertainty, called the metric coupling $\kappa = \{1/5, 2/5, 3/5, 4/5\}$ are shown. In each case the maximum uncertainty corresponds to $\frac{1}{Z} e_\kappa^{-1/2}(1)$.

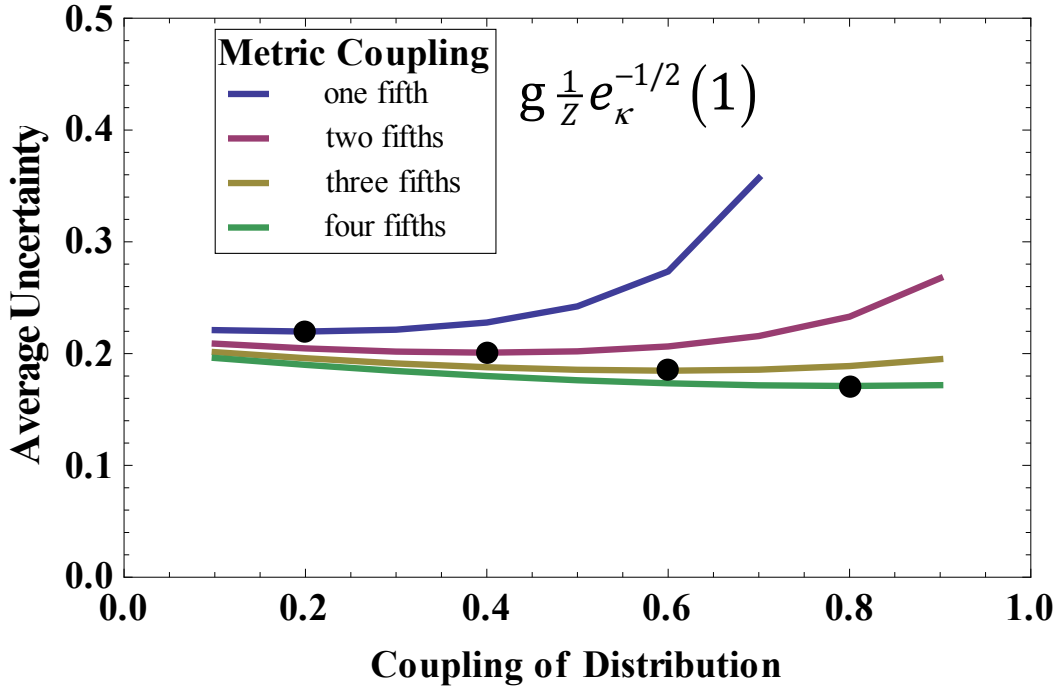

**Figure 2** The coupled average uncertainty as a function of the source of coupling is shown. The x-axis varies the coupling of a coupled Gaussian distribution with $\sigma = 1$. The four lines are for metric coupling of $\kappa = \{1/5, 2/5, 3/5, 4/5\}$, which translates into weighted generalized mean with $\frac{2\kappa}{1+\kappa}$. When the metric coupling equals the distribution coupling the average uncertainty is minimized and equals $\frac{1}{Z} e_\kappa^{-1/2}(1)$.

# 4 The coupled-entropy function

In defining the coupled log average in Equation (27) the role of the generalized entropy was implied but not explicitly expressed. In this section we show that the Renyi, Tsallis, and normalized Tsallis entropies can each be interpreted as the weighted generalized mean of a distribution translated to an entropy scale by generalized logarithms. Thus translated to the probability domain they are each identical expressions of average uncertainty. Use of the coupled logarithm leads to a new definition, called the coupled-entropy. Its properties with respect to the nonlinear coupling are compared with the Renyi and Tsallis generalizations.

In information theory, the negative logarithm of a probability is referred to as the surprisal. From Equation (16), the *coupled surprisal* (or generalized measure of information) is then

$$s_{\alpha,\kappa} \equiv -\ln_{-\alpha,\kappa} x \equiv \frac{1}{\alpha}\ln_\kappa x^{-\alpha} = \frac{1}{\alpha\kappa}\left(x^{\frac{-\alpha\kappa}{1+\kappa}} - 1\right), \text{ for } x > 0 \tag{37}$$

An alternative definition is possible based on the inverse of the coupled stretched exponential (20) taking $x/\sigma$ as the argument

$$\left(\ln_\kappa x^{-\alpha}\right)^{\frac{1}{\alpha}} \equiv \left(-\alpha \ln_{-\alpha,\kappa} x\right)^{\frac{1}{\alpha}} \equiv \left(\frac{1}{\kappa}\left(x^{\frac{-\alpha\kappa}{1+\kappa}} - 1\right)\right)^{\frac{1}{\alpha}}, \text{ for } x > 0. \tag{38}$$

While not investigated here, the alternative definition is of potential significance as an approach to the two-parameter generalized entropy function proposed by Hanel, et. al. [17,18] and related generalizations investigated in [37].

**Definition 4: Coupled Entropy** Given a discrete probability distribution $\mathbf{p} = \{p_i, i=1,...,N\}$ or a continuous density $f(x)$ and setting the Boltzmann constant to one, the *coupled entropy* is defined as

$$S_{\alpha,\kappa}(\mathbf{p}) \equiv -\ln_{-\alpha,\kappa}\left(\sum_{i=1}^{N} p_i^{1+\frac{\alpha\kappa}{1+\kappa}}\right)^{\frac{1+\kappa}{\alpha\kappa}} = -\sum_{i=1}^{N} P_i^{(\alpha,\kappa)} \ln_{-\alpha,\kappa} p_i, \tag{39}$$

$$S_{\alpha,\kappa}(f(x)) \equiv -\ln_{-\alpha,\kappa}\left(\int_{-\infty}^{\infty} f^{1+\frac{\alpha\kappa}{1+\kappa}}(x)dx\right)^{\frac{1+\kappa}{\alpha\kappa}} = -\int_{-\infty}^{\infty} f^{(\alpha,\kappa)}(x)\ln_{-\alpha,\kappa} f(x)dx. \tag{40}$$

The coupled-entropy expands to the following expression

$$S_{\alpha,\kappa}(\mathbf{p}) = -\ln_{-\alpha,\kappa}\left(\sum_{i=1}^{N} p_i^{1+\frac{\alpha\kappa}{1+\kappa}}\right)^{\frac{1+\kappa}{\alpha\kappa}}$$

$$= \frac{1}{\alpha\kappa}\left(\left(\sum_{i=1}^{N} p_i^{1+\frac{\alpha\kappa}{1+\kappa}}\right)^{\frac{1+\kappa-\alpha\kappa}{\alpha\kappa}\frac{1+\kappa}{1+\kappa}} - 1\right) \quad (41)$$

$$= \frac{1}{\alpha\kappa}\left(\left(\sum_{i=1}^{N} p_i^{1+\frac{\alpha\kappa}{1+\kappa}}\right)^{-1} - 1\right).$$

For the two principal cases of interest, namely for $\alpha = 1$ and $\alpha = 2$, the coupled-entropy becomes

$$S_{1,\kappa}(\mathbf{p}) = \frac{1}{\kappa}\left(\left(\sum_{i=1}^{N} p_i^{(1+2\kappa)/(1+\kappa)}\right)^{-1} - 1\right),$$

$$S_{2,\kappa}(\mathbf{p}) = \frac{1}{2\kappa}\left(\left(\sum_{i=1}^{N} p_i^{(1+3\kappa)/(1+\kappa)}\right)^{-1} - 1\right). \quad (42)$$

Both the Rényi and Tsallis entropies can be expressed in terms of the coupled average uncertainty (32). To facilitate the comparison between the entropy functions a substitution is made for the Rényi order and Tsallis index $q \to 1 + \frac{\alpha\kappa}{1+\kappa}$.

Table 1 summarizes the properties of the transformation from an average uncertainty to an entropy in terms of the moment of the average, and the power and normalization associated with the generalized logarithm,

$$S(\mathbf{p}, Norm, Power, Moment) = \frac{-1}{Norm}\left(GM^{Power} - 1\right), \quad GM = \left(\sum_{i=1}^{N} p_i^{Moment+1}\right)^{\frac{1}{Moment}}. \quad (43)$$

For the Rényi entropy ($R$) the relationship with the coupled average uncertainty is evident given its use of the natural logarithm

$$S_{\alpha,\kappa}^{R} = -\ln\left(\sum_{i=1}^{N} p_i^{1+\frac{\alpha\kappa}{1+\kappa}}\right)^{\frac{1+\kappa}{\alpha\kappa}}. \quad (44)$$

**Table 1 Comparison of Generalized Entropies** The generalized entropy functions are expressed in terms of three components, the moment of the generalized mean, the power and the normalization (norm) of the generalized logarithm.

| Entropy | Moment | Power | Norm. |
|---|---|---|---|
| Coupled | $\dfrac{\alpha\kappa}{1+\kappa}$ | $\dfrac{-\alpha\kappa}{1+\kappa}$ | $-\alpha\kappa$ |
| Tsallis | $\dfrac{\alpha\kappa}{1+\kappa}$ | $\dfrac{\alpha\kappa}{1+\kappa}$ | $\dfrac{\alpha\kappa}{1+\kappa}$ |
| Normalized Tsallis | $\dfrac{\alpha\kappa}{1+\kappa}$ | $\dfrac{-\alpha\kappa}{1+\kappa}$ | $\dfrac{-\alpha\kappa}{1+\kappa}$ |
| Renyi | $\dfrac{\alpha\kappa}{1+\kappa}$ | $\lim_{power \to 0}$ | $\lim_{norm \to 0}$ |
| Shannon | $\lim_{moment \to 0}$ | $\lim_{power \to 0}$ | $\lim_{norm \to 0}$ |

The coupled entropy is most closely related to the normalized Tsallis entropy (*NT*). Expressed in terms of the coupled surprisal (37) the normalized Tsallis entropy is seen to be the coupled entropy (39) multiplied by $(1+\kappa)$

$$S^{NT}_{\alpha,\kappa}(\mathbf{p}) \equiv \frac{-(1+\kappa)\sum_{i=1}^{N} p_i^{1+\frac{\alpha\kappa}{1+\kappa}} \ln_{-\alpha,\kappa} p_i}{\sum_{i=1}^{N} p_i^{1+\frac{\alpha\kappa}{1+\kappa}}} = -(1+\kappa)\ln_{-\alpha,\kappa}\left(\sum_{i=1}^{N} p_i^{1+\frac{\alpha\kappa}{1+\kappa}}\right)^{\frac{1+\kappa}{\alpha\kappa}}. \qquad (45)$$

Using the unnormalized form of Equation (45) for the Tsallis entropy definition, the lack of normalization is shown to be equivalent to using the opposite sign of the multiplier for the transformation

$$S^{T}_{\alpha,\kappa}(\mathbf{p}) \equiv -(1+\kappa)\sum_{i=1}^{N} p_i^{1+\frac{\alpha\kappa}{1+\kappa}} \ln_{-\alpha,\kappa} p_i = (1+\kappa)\ln_{\alpha,\kappa}\left(\sum_{i=1}^{N} p_i^{1+\frac{\alpha\kappa}{1+\kappa}}\right)^{\frac{1+\kappa}{\alpha\kappa}}. \qquad (46)$$

Derivation of the right-hand expression is shown in Appendix A. The coupled entropy (*C*), normalized Tsallis, and Tsallis entropy have then the following relationship

$$S^{C}_{\alpha,\kappa} = \frac{S^{NT}_{\alpha,\kappa}}{1+\kappa} = \frac{S^{T}_{\alpha,\kappa}}{(1+\kappa)\sum_{i=1}^{N} p_i^{1+\frac{\alpha\kappa}{1+\kappa}}}. \qquad (47)$$

Figure 3 and Figure 4 show comparisons of the entropy functions for the coupled exponential distribution and the coupled Gaussian distribution respectively. The computations for the two figures are summarized in Appendix B. For the coupled exponential distribution (Equation (21) with $\alpha=1$) with the scale $\sigma$ equal to one, both the coupled and Tsallis entropy are equal to one. In this case, the normalized Tsallis entropy is linear and similar to the Shannon entropy. As shown on the right of Figure 3, even as the scale varies the coupled entropy does not vary significantly with the coupling. This is suggestive that the coupled entropy represents a stabilization of the entropy as a function of the coupling, in contrast to the Shannon and Renyi entropies which have close to linear growth with the coupling. However, this property does not hold for the coupled Gaussian distribution (Equation (21) with $\alpha=2$) shown in Figure 4. Instead, the coupled entropy is quite sensitive to changes in both the coupling and the scale. With scale equal to one, the coupled entropy is similar to the Shannon entropy, particularly when $0<\kappa<0.5$. For $\kappa>0.5$, the difference is more evident as the coupled entropy is slightly super-linear and the Shannon entropy is slightly sub-linear. For this case, the normalized Tsallis entropy is close to quadratic with the coupling. The Renyi entropy has close to linear growth for both distributions, but with a slope that is less than the Shannon entropy. The sensitivity of the coupled entropy with $\alpha\neq1$ would be reduced if Equation (38) is used for the generalized logarithm, since the root with respect to $\alpha$ is taken.

# 5 Conclusion

The nonlinear statistical coupling of states is defined in terms of a source of nonlinear coupling $\kappa$ in the system. Positive coupling models non-stationary systems in which the variance fluctuates. The negative domain of $\kappa$, which has decreased variability, is limited to $(-1,0)$; and $\kappa=0$ is the linear domain with exponential distributions. The maximum entropy distribution for these systems given a constrained location and scale is the coupled-exponential distributions, defined in Equation (21). In addition to the nonlinear coupling these distributions depend on the stability index of the stretched exponential $e^{-x^\alpha/\alpha}$. The coupled sum is defined by the nonlinear coupling, while the coupled product is a function of the coupling, the stability index and the dimensions $\dfrac{-\alpha\kappa}{1+d\kappa}$. In this paper the one-dimensional case $d=1$ was examined.

In this paper the coupled average uncertainty in the probability/density space is introduced as the weighted generalized mean with a power equal to $\frac{\alpha\kappa}{1+\kappa}$. The expression is derived from and is expressible in terms of deformation of the coupled product (18) and the coupled power (19) functions

$$P_{\kappa avg} \equiv \prod_{i=1}^{N} {}_{\otimes_{\alpha\kappa_1}} p_i^{\otimes_{\kappa}^{p_i^{(-\alpha\kappa_1)}}} = \left(\sum_{i=q}^{N} p_i^{1+\frac{\alpha\kappa}{1+\kappa}}\right)^{\frac{1+\kappa}{\alpha\kappa}}. \qquad (48)$$

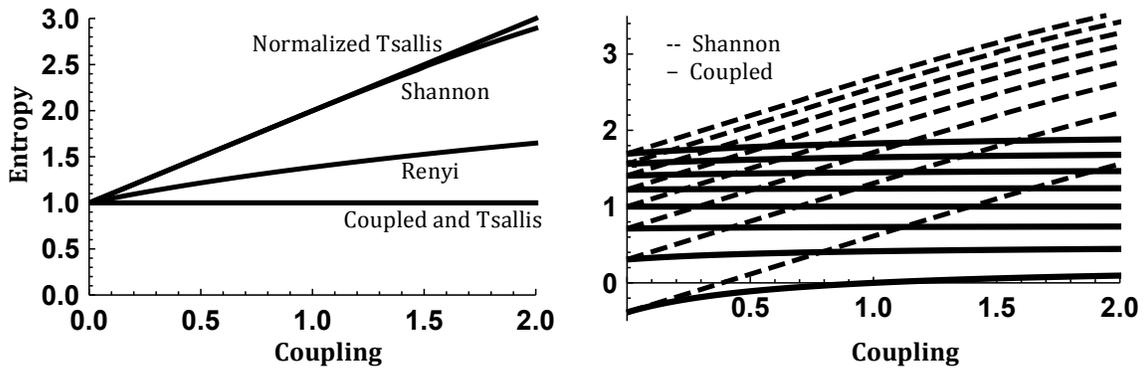

**Figure 3 Comparison of Entropy Functions for Coupled Exponential Distribution** Left) The entropy versus coupling with the scale equal to one. The Coupled and Tsallis entropies are both equal to 1. The Shannon and Renyi entropy have close to linear growth. The Normalized Tsallis entropy is similar to Shannon entropy. Right) Comparison of Shannon (dashed) and Coupled (solid) entropy for $\sigma = \{0.25, 0.5, 0.75, 1, 1.25, 1.5, 1.75, 2\}$. Both entropies increase with the scale $\sigma$. The Coupled entropy does not vary significantly with the coupling, while the Shannon entropy has close to linear growth.

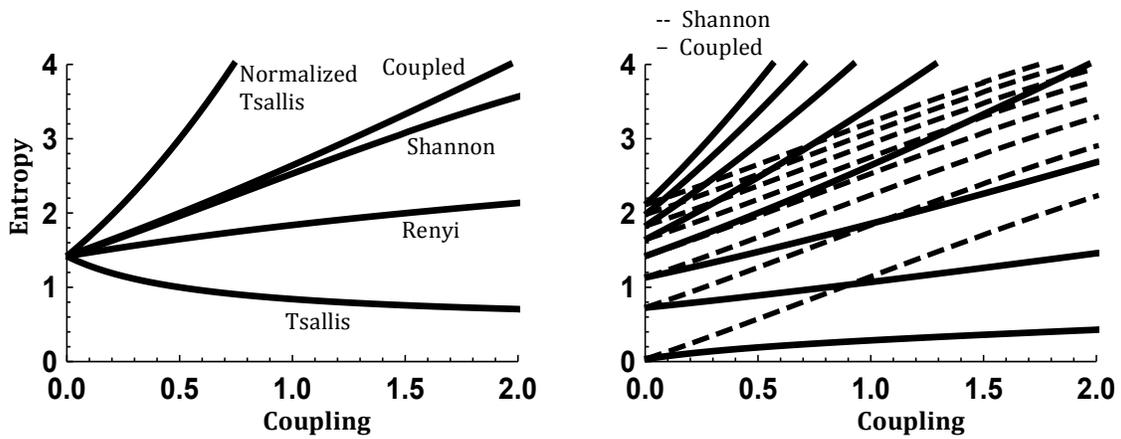

**Figure 4 Comparison of Entropy Functions for Coupled Gaussian Distribution** Left) The entropy versus coupling with the scale equal to one. The Coupled and Shannon entropy are similar with close to linear growth. The Renyi entropy also has close linear growth but with a smaller slope. The Normalized Tsallis entropy has grow close to the square of the coupling, while the Tsallis entropy decays. Right) Comparison of Shannon (dashed) and Coupled (solid) entropy for $\sigma = \{0.25, 0.5, 0.75, 1, 1.25, 1.5, 1.75, 2\}$. Both entropies increase with the scale $\sigma$. The Coupled entropy has increased sensitivity with $\sigma$ as the coupling increase, while the Shannon entropy has close to linear growth.

For $\kappa = 0$ this is the weighted geometric mean in which each state is treated as independent and thus multiplied, and the dependency from the probabilities summing to one is contained in the weight of each probability. When $\kappa \neq 0$ the coupled product models the nonlinear dependency between the probabilities of the distribution and the weights are now the deformed coupled probability. Importantly, for the coupled exponential distributions the coupled average probability is equal to the density at the location plus the scale $P_{\kappa avg} = \frac{1}{Z} e_\kappa^{-1/\alpha}(1)$. The coupled average probability has been successfully utilized in the analysis of information fusion algorithms [8,41] and is recommended as a simple approach to understanding the average uncertainty for nonlinear systems.

Upon transformation to the probability domain the Tsallis, normalized Tsallis and Rényi entropies each fulfill the relationship between the coupled average uncertainty and the density at the width of the coupled distributions. The potentially improved properties of the coupled-entropy are achieved by requiring the coupled logarithm with $\alpha = 1$ to have a unit integral. The coupled entropy is defined as the transformation of $P_{\kappa avg}$ to the entropy scale by the coupled logarithm,

$$S_{\alpha,\kappa}(\mathbf{p}) \equiv -\ln_{-\alpha,\kappa} \left( \sum_{i=1}^{N} p_i^{1+\frac{\alpha\kappa}{1+\kappa}} \right)^{\frac{1+\kappa}{-\alpha\kappa}} = \frac{1}{\alpha\kappa} \left( \left( \sum_{i=1}^{N} p_i^{(1+\alpha\kappa)/(1+\kappa)} \right)^{-1} - 1 \right). \tag{49}$$

For the coupled exponential distribution, the matching coupled entropy does not vary significantly with the degree of coupling. In contrast, for the Gaussian distribution, the matching coupled entropy is sensitive to both the coupling and scale of the distribution. Further analysis and experimentation is recommended to determine how the coupled entropy can be used to model the effects of non-stationary uncertainty in complex thermodynamic and information systems.

# Appendices

## A. Derivation of coupled entropy from average uncertainty

The coupled algebra can be used to derive the equivalency between the two expressions for the coupled entropy

$$S_{\alpha,\kappa}(\mathbf{p}) \equiv -\ln_{-\alpha,\kappa}\left(\sum_{i=1}^{N} p_i^{1+\frac{\alpha\kappa}{1+\kappa}}\right)^{\frac{1+\kappa}{\alpha\kappa}} = -\sum_{i=1}^{N} P_i^{(\alpha,\kappa)} \ln_{-\alpha,\kappa} p_i$$

where $P_i^{(\alpha,\kappa)}$ is the coupled probability (22)

$$P_i^{(\alpha,\kappa)} = \frac{p_i^{1+\frac{\alpha\kappa}{1+\kappa}}}{\sum_{j=1}^{N} p_j^{1+\frac{\alpha\kappa}{1+\kappa}}}$$

Starting with the expression on the right, the coupled power is used to move the coupled probability from a multiplier of the generalized logarithm to a power of the argument using the property $b \ln_{-\alpha,\kappa} x = \ln_{-\alpha,\kappa} x^{\otimes^b_{-\alpha\kappa/(1+\kappa)}}$

$$S_{\alpha,\kappa}(\mathbf{p}) = -\sum_{i=1}^{N} \ln_{-\alpha,\kappa} p_i^{\otimes^{P_i^{(-\alpha\kappa/(1+\kappa))}}_{-\alpha\kappa/(1+\kappa)}}.$$

Expanding the coupled power and using the property $\sum_{i=1}^{N} \ln_{-\alpha,\kappa} x = \ln_{-\alpha,\kappa} \prod_{i=1}^{N} \otimes_{\frac{-\alpha\kappa}{1+\kappa}} x$,

$$S_{\alpha,\kappa}(\mathbf{p}) = -\ln_{-\alpha,\kappa} \prod_{i=1}^{N} \otimes_{\frac{-\alpha\kappa}{1+\kappa}} \left(P_i^{(\alpha,\kappa)} p_i^{\frac{-\alpha\kappa}{1+\kappa}} - \left(P_i^{(\alpha,\kappa)} - 1\right)\right)^{\frac{1+\kappa}{-\alpha\kappa}}.$$

Expanding the coupled product and reducing $\sum_{i=1}^{N}\left(P_i^{(\alpha,\kappa)} - 1\right) = 1 - N$,

$$S_{\alpha,\kappa}(\mathbf{p}) = -\ln_{-\alpha,\kappa}\left(\sum_{i=1}^{N} P_i^{(\alpha,\kappa)} p_i^{\frac{-\alpha\kappa}{1+\kappa}}\right)^{\frac{1+\kappa}{-\alpha\kappa}} = -\ln_{-\alpha,\kappa}\left(\frac{\sum_{i=1}^{N} p_i}{\sum_{i=1}^{N} p_i^{1+\frac{\alpha\kappa}{1+\kappa}}}\right)^{\frac{1+\kappa}{-\alpha\kappa}}$$

$$= -\ln_{-\alpha,\kappa}\left(\sum_{i=1}^{N} p_i^{1+\frac{\alpha\kappa}{1+\kappa}}\right)^{\frac{1+\kappa}{\alpha\kappa}}.$$

Since the normalized Tsallis entropy is just a multiple of the coupled entropy, the derivation is not repeated. The derivation of the Tsallis entropy starts with the middle expression of equation (46)

$$S^T_{\alpha,\kappa}(\mathbf{p}) \equiv -(1+\kappa)\sum_{i=1}^{N} p_i^{1+\frac{\alpha\kappa}{1+\kappa}} \ln_{-\alpha,\kappa} p_i.$$

The relationship is expressed in terms of the coupled product and coupled power

$$S^T_{\alpha,\kappa}(\mathbf{p}) = -(1+\kappa)\ln_{-\alpha,\kappa} \prod_{i=1}^{N} {}_{\otimes_{\frac{-\alpha\kappa}{1+\kappa}}} p_i^{\otimes p_i^{1+\frac{\alpha\kappa}{1+\kappa}}}_{-\alpha\kappa/(1+\kappa)}.$$

Expanding the coupled power and the coupled power

$$S^T_{\alpha,\kappa}(\mathbf{p}) = -(1+\kappa)\ln_{-\alpha,\kappa}\left(\sum_{i=1}^{N}\left(p_i - \left(p_i^{1+\frac{\alpha\kappa}{1+\kappa}} - 1\right)\right) - (N-1)\right)^{\frac{1+\kappa}{-\alpha\kappa}}.$$

Simplifying and using the relationship $-\ln_{-\alpha,\kappa} x^{-b} = \ln_{\alpha,\kappa} x^b$

$$S^T_{\alpha,\kappa}(\mathbf{p}) = (1+\kappa)\ln_{\alpha,\kappa}\left(\sum_{i=1}^{N} p_i^{1+\frac{\alpha\kappa}{1+\kappa}}\right)^{\frac{1+\kappa}{\alpha\kappa}}.$$

## B. Comparison of generalized entropy functions

The comparison of entropy functions for the Coupled Exponential and Coupled Gaussian distributions shown in Figure 3 and Figure 4 includes the following computations. The distributions analyzed are

$$expdist_\kappa(x;\mu,\sigma) \equiv \tfrac{1}{\sigma}\left(1+\kappa\tfrac{x-\mu}{\sigma}\right)^{-\left(\frac{1}{\kappa}+1\right)} \quad (50)$$

$$G_\kappa(x;\mu,\sigma) \equiv \frac{\sqrt{\kappa}\,\Gamma\!\left[\frac{1+\kappa}{2\kappa}\right]}{\sqrt{\pi}\,\Gamma\!\left[\frac{1}{2\kappa}\right]}\left(1+\kappa\left(\frac{x-\mu}{\sigma}\right)^2\right)^{-\frac{1}{2}\left(\frac{1}{\kappa}+1\right)}. \quad (51)$$

The computations were completed using *Mathematica*™. For the Coupled, Tsallis, and Normalized Tsallis entropy with matching coupling, closed form solutions exist. The entropy functions of the coupled exponential with $\sigma = 1$ reduce to a linear expression.

a) Coupled Entropy

$$S^C_{1,\kappa}(expdist_\kappa(x;\mu,\sigma)) = \frac{-1+(1+\kappa)\sigma^{\frac{\kappa}{1+\kappa}}}{\kappa} \quad (52)$$

$$S^C_{1,\kappa}(expdist_\kappa(x;\mu,1)) = 1$$

$$S_{2,\kappa}^{C}\left(G_{\kappa}(x;\mu,\sigma)\right)=\frac{-\kappa+(1+\kappa)\left(\kappa\pi^{\kappa}\right)^{\frac{1}{1+\kappa}}\left(\sigma\frac{\Gamma\left(\frac{1}{2\kappa}\right)}{\Gamma\left(\frac{1+\kappa}{2\kappa}\right)}\right)^{\frac{2\kappa}{1+\kappa}}}{2\kappa^{2}} \quad (53)$$

b) Tsallis Entropy

$$S_{1,\kappa}^{T}\left(expdist_{\kappa}(x;\mu,\sigma)\right)=\frac{1+\kappa-\sigma^{\frac{-\kappa}{1+\kappa}}}{\kappa} \quad (54)$$

$$S_{1,\kappa}^{T}\left(expdist_{\kappa}(x;\mu,1)\right)=1$$

$$S_{2,\kappa}^{T}\left(G_{\kappa}(x;\mu,\sigma)\right)=\frac{1}{2}\left(1+\frac{1}{\kappa}-\left(\kappa\pi^{\kappa}\right)^{\frac{-1}{1+\kappa}}\left(\frac{\sigma\Gamma\left[\frac{1}{2\kappa}\right]}{\Gamma\left[\frac{1+\kappa}{2\kappa}\right]}\right)^{\frac{-2\kappa}{1+\kappa}}\right) \quad (55)$$

c) Normalized Tsallis

$$S_{1,\kappa}^{NT}\left(expdist_{\kappa}(x;\mu,\sigma)\right)=\frac{(1+\kappa)\left(-1+(1+\kappa)\sigma^{\frac{\kappa}{1+\kappa}}\right)}{\kappa} \quad (56)$$

$$S_{1,\kappa}^{NT}\left(expdist_{\kappa}(x;\mu,1)\right)=1+\kappa$$

$$S_{2,\kappa}^{NT}\left(G_{\kappa}(x;\mu,\sigma)\right)=\frac{(1+\kappa)\left(-\kappa+(1+\kappa)\left(\kappa\pi^{\kappa}\right)^{\frac{1}{1+\kappa}}\left(\frac{\sigma\Gamma\left(\frac{1}{2\kappa}\right)}{\Gamma\left(\frac{1+\kappa}{2\kappa}\right)}\right)^{\frac{2\kappa}{1+\kappa}}\right)}{2\kappa^{2}} \quad (57)$$

Closed form solutions do not exist for the Shannon and Renyi entropies. The numerical integration of the Shannon entropy has potential errors for a) small values of $\kappa$ when the integration extends to regions where $-\ln f(x)^{f(x)}$ is undefined or b) large values of $\kappa$ when the integration is truncated before $-\ln f(x)^{f(x)}$ has decayed sufficiently to be excluded. For this reason the Shannon entropy was determined as a piecewise numerical integration with the following limits

$$\begin{array}{ll} x \text{ Limit} & \text{Coupling} \\ 100 & 0.00<\kappa<0.09 \\ 1000 & 0.09\leq\kappa<0.74 \\ 10{,}000 & 0.74\leq\kappa<1.50 \\ 15{,}000 & 1.50\leq\kappa\leq 2.00. \end{array} \quad (58)$$

# References

[1] T.M. Cover, J.A. Thomas, Elements of information theory, John Wiley & Sons, 2012.

[2] A.I. Khinchin, Mathematical Foundations of Statistical Mechanics, Courier Corporation, 1949.

[3] A.I. Khinchin, Mathematical Foundations of Information Theory, Courier Corporation, 1957.

[4] A. Rényi, On measures of entropy and information, Fourth Berkeley Symp. Math. Stat. Probab. 1 (1961) 547–561.

[5] C. Tsallis, Possible generalization of Boltzmann-Gibbs statistics, J. Stat. Phys. 52 (1988) 479–487.

[6] S.I. Amari, H. Nagaoka, Methods of Information Geometry, Oxford University Press, 2000.

[7] C. Anteneodo, C. Tsallis, A.S. Martinez, Risk aversion in trade transactions, Eur. Lett. 59 (2002) 635–641.

[8] K.P. Nelson, B.J. Scannell, H. Landau, A Risk Profile for Information Fusion Algorithms, Entropy. 13 (2011) 1518–1532.

[9] I. Santamaria, D. Erdogmus, J.C. Principe, Entropy minimization for supervised digital communications channel equalization, IEEE Trans. Signal Process. 50 (2002) 1184–1192. doi:10.1109/78.995074.

[10] W.R. Carson, M. Chen, M.R.D. Rodrigues, R. Calderbank, L. Carin, Communications-Inspired Projection Design with Application to Compressive Sensing, SIAM J. Imaging Sci. 5 (2012) 1185–1212. doi:10.1137/120878380.

[11] E.P. Borges, C. Tsallis, G.F.J. Ananos, P.M.C. de Oliveira, Nonequilibrium probabilistic dynamics of the logistic map at the edge of chaos, Phys. Rev. Lett. 89 (2002) 254103.

[12] G. Livadiotis, D.J. McComas, Beyond kappa distributions: Exploiting Tsallis statistical mechanics in space plasmas, J. Geophys. Res. 114 (2009) A11105. doi:10.1029/2009JA014352.

[13] L.F. Burlaga, A.F. Vinas, Triangle for the entropic index q of non-extensive statistical mechanics observed by Voyager 1 in the distant heliosphere, Phys. A Stat. Mech. Its Appl. 356 (2005) 375–384.

[14] S. Chatrchyan, V. Khachatryan, A.M. Sirunyan, A. Tumasyan, W. Adam, T. Bergauer, et al., Charged particle transverse momentum spectra in pp collisions at $\sqrt{s}=0.9$ and 7 TeV, J. High Energy Phys. 2011 (2011) 86. doi:10.1007/JHEP08(2011)086.

[15] V. Khachatryan, A.M. Sirunyan, A. Tumasyan, W. Adam, T. Bergauer, M. Dragicevic, et al., Transverse-momentum and pseudorapidity distributions of charged hadrons in pp collisions at square root of s = 7 TeV., Phys. Rev. Lett. 105 (2010) 022002. doi:10.1103/PhysRevLett.105.022002.

[16] N. Canosa, R. Rossignoli, Generalized nonadditive entropies and quantum entanglement, Phys. Rev. Lett. 88 (2002) 170401.

[17] R. Hanel, S. Thurner, M. Gell-Mann, Generalized entropies and the transformation group of superstatistics, Proc. Natl. Acad. Sci. 108 (2011) 6390.

[18] R. Hanel, S. Thurner, A comprehensive classification of complex statistical systems and an axiomatic derivation of their entropy and distribution functions, EPL (Europhysics Lett. 93 (2011) 20006.

[19] A. Rényi, On the Foundations of Information Theory, Rev. Int. Stat. Inst. 33 (1965) 1–14.

[20] C. Tsallis, Introduction to Nonextensive Statistical Mechanics: Approaching a Complex World, Springer Verlag, 2009.

[21] K.P. Nelson, S. Umarov, Nonlinear statistical coupling, Phys. A Stat. Mech. Its Appl. 389 (2010) 2157–2163.

[22] T.J. Arruda, R.S. González, C.A.S. Terçariol, A.S. Martinez, Arithmetical and geometrical means of generalized logarithmic and exponential functions: generalized sum and product operators, Phys. Lett. A. 372 (2008) 2578–2582.

[23] Q.A. Wang, L. Nivanen, A. Le Mehaute, M. Pezeril, On the generalized entropy pseudoadditivity for complex systems, J. Phys. A. Math. Gen. 35 (2002) 7003–7007.

[24] G. Kaniadakis, M. Lissia, A.M. Scarfone, Two-parameter deformations of logarithm, exponential, and entropy: A consistent framework for generalized statistical mechanics, Phys. Rev. E. 71 (2005) 046128. doi:10.1103/PhysRevE.71.046128.

[25] K.S. Narendra, K. Parthasarathy, Identification and control of dynamical systems using neural networks., IEEE Trans. Neural Netw. 1 (1990) 4–27. doi:10.1109/72.80202.

[26] S.H. Strogatz, Exploring complex networks., Nature. 410 (2001) 268–76. doi:10.1038/35065725.

[27] T. Oikonomou, Tsallis, Renyi and nonextensive Gaussian entropy derived from the respective multinomial coefficients, Phys. A Stat. Mech. Its Appl. 386 (2007) 119–134.

[28] K.P. Nelson, A definition of the coupled-product for multivariate coupled-exponentials, Phys. A Stat. Mech. Its Appl. 422 (2015) 187–192. doi:10.1016/j.physa.2014.12.023.

[29] E.P. Borges, A possible deformed algebra and calculus inspired in nonextensive thermostatistics, Phys. A Stat. Mech. Its Appl. 340 (2004) 95–101.

[30] C. Anteneodo, C. Tsallis, Multiplicative noise: A mechanism leading to nonextensive statistical mechanics, J. Math. Phys. 44 (2003) 5194.

[31] C. Beck, E.G.D. Cohen, Superstatistics, Phys. A Stat. Mech. Its Appl. 322 (2003) 267–275.

[32] G. Wilk, Z. Włodarczyk, Interpretation of the Nonextensivity Parameter q in Some Applications of Tsallis Statistics and Lévy Distributions, Phys. Rev. Lett. 84 (2000) 2770–2773. doi:10.1103/PhysRevLett.84.2770.

[33] W. Gossett, The Application of the "Law of Error" to the Work of the Brewery, Dublin, 1904.

[34] A. Shah, A. Wilson, Z. Ghahramani, Student-t processes as alternatives to Gaussian processes, arXiv Prepr. arXiv1402.4306. (2014).

[35] A.C. de Souza, C. Tsallis, Student's t- and r-distributions: Unified derivation from an entropic variational principle, Phys. A Stat. Mech. Its Appl. 236 (1997) 52–57. doi:10.1016/S0378-4371(96)00395-0.

[36] S. Umarov, C. Tsallis, On multivariate generalizations of the q-central limit theorem consistent with nonextensive statistical mechanics, AIP Conf. Proc. 965 (2007) 34.

[37] G. Kaniadakis, M. Lissia, A.M. Scarfone, Two-parameter deformations of logarithm, exponential, and entropy: A consistent framework for generalized statistical mechanics, Phys. Rev. E. 71 (2005) 46128.

[38] S. Umarov, C. Tsallis, M. Gell-Mann, S. Steinberg, Generalization of symmetric alpha-stable Lévy distributions for q>1., J. Math. Phys. 51 (2010) 33502. doi:10.1063/1.3305292.

[39] C. Tsallis, A.R. Plastino, R.F. Alvarez-Estrada, Escort mean values and the characterization of power-law-decaying probability densities, J. Math. Phys. 50 (2009) 43303.

[40] C. Lassner, R. Lienhart, Norm-Induced Entropies for Decision Forests, in: 2015 IEEE Winter Conf. Appl. Comput. Vis., IEEE, 2015: pp. 968–975. doi:10.1109/WACV.2015.134.

[41] K.P. Nelson, Reduced perplexity: Uncertainty measures without entropy, in: Recent Innov. Info-Metrics, Washington, D.C., 2014.